\documentclass[11pt]{article}
\usepackage{amsmath,amssymb}

\newcommand{\Tbar}{\ensuremath{\overline{\mathcal{T}}}}

\newcommand{\caT}{\mathcal T}

\begin{document}

\title{Families of Riemann surfaces and Weil-Petersson geometry}         
\author{Scott A. Wolpert}        
\date{April 11, 2008}          
\maketitle

\section{Considerations}

The investigation of the geometry and deformation theory of Riemann surfaces is long recognized for involving a wide range of concepts and a wide range of techniques.  For surfaces of negative Euler characteristic, the Weil-Petersson metric provides a Hermitian structure for the space of infinitesimal variations of hyperbolic metrics.  The Weil-Petersson geometry is recognized for providing important information on the geometry of the Teichm\"{u}ller space and the moduli space of Riemann surfaces.  In less than a decade there have been a collection of breakthrough works on Weil-Petersson geometry.  The works are based on a range of approaches and so involve matters of interest and further study for a broad audience.

The plan for the lecture topics is presented in Section \ref{lect}.  The approach will include the fundamentals and central underlying arguments for three areas: explorations of the $CAT(0)$ geometry, applications of volume integrals and applications of curvature.  Supporting goals for the approach are to present results that are based on standardized techniques and to present results that lead directly into research questions.   

\section{Weil-Petersson background and recent results}
By the Uniformization Theorem a genus $g$ with $n$ punctures, $2g-2+n>0$,  Riemann surface $R$ has a unique hyperbolic metric $ds^2$.  By Riemann-Roch the associated space of holomorphic quadratic differentials $Q(R)$ on $R$ with at most simple poles at punctures has dimension $3g-3+n$.  By Kodaira-Spencer deformation theory and Serre duality $Q(R)$ is the dual of the space of infinitesimal deformations of $R$. Weil introduced a Hermitian pairing for the space of infinitesimal deformations by defining
\[
\langle\varphi,\psi\rangle = \int_R \varphi\overline{\psi}(ds^2)^{-1}\mbox{ for }\varphi,\psi\in Q(R),\quad \cite{Ahsome}.
\]
The Teichm\"{u}ller space $\caT$ is the space of homotopy-marked hyperbolic metrics for a surface of genus $g$ with $n$ punctures.  The Weil-Petersson (WP) metric for $\caT$ is invariant under the action of the mapping class group $MCG$ (the mapping class group is 
$\operatorname{Homeo}^+(R)/\operatorname{Homeo}_0(R)$). The metric is K\"{a}hler, non complete with negative sectional curvatures $\kappa$ with $\sup_{\caT}\kappa=0$ (except for $\dim_{\mathbb C} \caT =1$) and $\inf_{\caT}\kappa=-\infty$.  The metric completion of $\caT$ is the augmented Teichm\"{u}ller space $\Tbar$, the space of homotopy-marked noded Riemann surfaces.  $\Tbar$ is a disjoint union of $\caT$ and lower dimensional Teichm\"{u}ller spaces.  With the completion metric $\Tbar$ is an example of a $CAT(0)$ metric space, a simply connected generalized space of negative curvature, \cite{DW2, Wlcomp}.  Alexandrov and others developed a collection of techniques for understanding $CAT(0)$ spaces, \cite{BH}.  The topological pants graph $P(R)$ has a vertex for each collection of $3g-3+n$ distinct free homotopy classes of non trivial, non peripheral, mutually disjoint simple closed curves.  Vertices of $P(R)$ are connected by an edge provided the corresponding sets of free homotopy classes differ by replacing a single curve with a once or twice intersecting curve (a Hatcher-Thurston move).  Brock discovered that Teichm\"{u}ller space with the WP metric is quasi isometric to $P(R)$ with the unit-edge metric, \cite{Brkwp}.  At large scale on Teichm\"{u}ller space WP distance is combinatorially determined.

The importance of the Weil-Petersson metric begins with the natural family of Riemann surfaces over Teichm\"{u}ller space being the universal family of homotopy-marked hyperbolic metrics.  There are wide ranging applications of the geometry.  The WP K\"{a}hler geometry provides a complex differential geometric proof with the Kodaira Embedding Theorem that the Deligne-Mumford compactification $\overline{\caT/MCG}$ is projective algebraic, \cite{Wlpos}.  The metric convexity provides for Nielsen Realization: finite subgroups of MCG acting on $\caT$ have non empty fixed-point sets, \cite{Wlnielsen}.  Mirzakhani's WP volume recursion (based on symplectic reduction) in topological type $g,n$ provides for a solution of the Witten-Kontsevich conjecture,\cite{Mirwitt, Mirvol}. Li-Sun-Yau use the WP Ricci form as a beginning reference metric for constructing and estimating the K\"{a}hler-Einstein metric on $\caT$, \cite{LSY1,LSY2}.  They establish comparability of the K\"{a}hler-Einstein metric to the five classical complete metrics (Kobayashi-Teichm\"{u}ller, Carath\'{e}odory, asymptotic Poincar\'{e}, McMullen and Bergman). Up to comparability the six metrics prospectively share properties.  

There are also wide ranging descriptions of the WP metric.  For a deformation family of Riemann surfaces  $R_t$ beginning with $R$, consider the energy $E(f_t)$ of the harmonic map $f_t:R\rightarrow R_t$.  Up to a universal constant the Hessian of $E(f_t)$ is the WP metric at $R$, \cite{Wfharm}.  For the quasiFuchsian group $\Gamma_t$ representing the pair of Riemann surfaces $(R,R_t)$ then up to a constant the Hessian of the Hausdorff dimension of the limit set of $\Gamma_t$ is also the WP metric, \cite{McMthr}.  There is also a characterization in terms of geodesic-length functions.  For a geodesic $\gamma$ on $R$  let $\ell_{\gamma}(R_t)$ be the length of the unique hyperbolic metric geodesic in the corresponding free homotopy class on $R_t$.  Let $\{\gamma_j\}$ be a sequence of geodesics with length-normalized lifts to the unit-tangent bundle $TR$ tending weak$^*$ to the uniform distribution. The $\ell_{\gamma_j}$ Hessians tend to the WP metric as follows (for $n=0$; Thurston's random geodesic metric is the WP metric)
\[
\langle\ ,\ \rangle\,=\, 3\pi(g-1)\lim_j \frac{1}{\ell_j}\operatorname{Hess}\ell_{\gamma_j}\quad\cite{Wlthur}.
\]

\section{Weil-Petersson research themes}
The deformation theory of Riemann surfaces and hyperbolic metrics provides an entree to deformation theory of higher dimensional manifolds.  Deformations of Riemann surfaces represent a basic case of the application of: Kodaira-Spencer theory; the prescribed curvature equation; geometric structures described by developing maps; $\operatorname{Hom}(\pi_1(R),G)/G$ ($G=PSL(2;\mathbb R)$ for hyperbolic metrics) and for Quillen's metric approach to the local index theorem for the complex differential (Cauchy-Riemann operator) $\overline{\partial}$.  

A theme of investigation is intersection theory and understanding curvature.  Intersection theory of cycles on the Deligne-Mumford compactification $\overline{\caT/MCG}$ was employed, \cite{Wlhomol}, to evaluate the WP class and is important to Penner's resolution of Witten's conjecture on an explicit cycle Poincar\'{e} dual to the WP K\"{a}hler form, \cite{Pen2frm}.   The sectional curvature was originally considered by Royden, Tromba and the author, \cite{Royicm, Trcurv, Wlchern}.  Recent results include that the WP metric is Gromov hyperbolic iff $\dim_{\mathbb C}\caT \le 2$, \cite{BF}.  Zheng Huang has been effecting a detailed analysis of WP curvature, \cite{Zh, Zh2, Zh3, Zh4}.  While the WP volume grows exponentially in $g$, he finds that on the thick part of the moduli spaces the curvatures are bounded.  A precise application of curvature and the hyperbolic metric is given Freixas' development of arithmetic intersection theory  following Gillet-Soul\'{e} and Bismut for Riemann surfaces with punctures, \cite{GFM}.  A consequence is the exact formula for the Selberg zeta function
\[
Z'(1)\,=\,4\pi^{5/3}\Gamma_2(1/2)^{-8/3}
\]
for the principal subgroup $\Gamma(2)\subset SL(2;\mathbb Z)$ and the Barnes double gamma function.  The works of Takhtajan and Zograf present an original application of Quillen's metric for calculating the  Chern form of the determinant line bundle, \cite{TZcomp, TZ, TZpunc}.  The determinant is given by a special value of the Selberg zeta function.    

An important theme of investigation is the WP symplectic geometry. An original result is the duality $\omega(t_{\gamma},\ )\,=\,-\frac12 d\ell_{\gamma}$ of Fenchel-Nielsen infinitesimal deformations and geodesic-length functions relative to the WP symplectic form $\omega$, \cite{Wlsymp}.  The symplectic geometry is manifested in the formula for Fenchel-Nielsen twist-length coordinates
\[
\omega\,=\, \frac12 \sum_{\gamma_j}d\ell_j\wedge d \tau_j, \quad\cite{Wldtau},
\]
where $\{\gamma_j\}$ gives a pants decomposition (the expression is universal even though there are an infinite number of topologically distinct pants decompositions).   Goldman more generally considered trace-functions on the representation spaces $\operatorname{Hom}(\pi_1(R),G)/G$ as Hamiltonians relative to the cup product for the group-cohomology description of the tangent space, \cite{Gdsymp}.   

In a series of innovative papers Mirzakhani presented calculations of WP integrals exploiting the underlying symplectic geometry, \cite{Mirvol, Mirgrow, Mirwitt}.    She also established a {\em prime geodesic theorem} for simple geodesics $\#\{\gamma\mid\ell_{\gamma}(R)\le L\}\,\sim\,c_R L^{6g-6+2n}$, for $L$ large, and a proof of the Witten Kontsevich formulas for intersections of tautological classes on the Deligne-Mumford compactification $\overline{\caT/MCG}$.  
Her work is an entree to several matters.  She combines a generalization of McShane's universal length-sum identity, recursions of integrals, Thurston's symplectic geometry for the space of measured geodesic laminations and symplectic reduction of quasi-free $S^1$ actions to obtain a collection of new results.  She finds volume formulas for the moduli spaces of hyperbolic surfaces with geodesic boundary including the formula
\[
vol_{WP}\,=\,\frac{\pi^2}{6}\,+\,\frac{b^2}{24}
\]
for the moduli space of tori with length $b$ boundary.

Another theme of WP investigation is the metric space geometry and the behavior of geodesics.  As above the augmented Teichm\"{u}ller space $\Tbar$ is a $CAT(0)$ metric space with $\caT,\,\Tbar$ quasi isometric to the pants graph $P(R)$.  Recently the Alexandrov tangent cones of $\Tbar$ have been described, \cite{Wlbhv}.  The description provides for applications and for considering combinatorial harmonic maps into $\Tbar$.  A refined study of WP geodesics is now underway.   Brock-Masur-Minsky have introduced a notion of {\em ending lamination} for infinite WP rays and apply their notion to embed asymptote classes (of recurrence rays on $\caT/MCG$) into the Gromov boundary of the curve complex $C(R)$ (the complex with $k$ simplices being collections of $k+1$ distinct free homotopy classes of non trivial, non peripheral, mutually disjoint simple closed curves on $R$), \cite{BMM}.  The authors apply their analysis to show that the WP geodesic flow is topologically transitive.  Pollicott-Weiss-Wolpert use a form of the classical closing lemma for geodesics on surfaces to obtain a geometric elementary proof for $\dim_{\mathbb C}\caT=1$, \cite{PWWl}.    

\section{Weil-Petersson research techniques}
A range of techniques is applied for the study of variations of hyperbolic metrics and WP geometry.   A majority of approaches have a basis in complex analysis and complex differential geometry.  Deformation theory plays the largest role.

An application of the Dolbeault complex is that the potential equation $\overline{\partial}f=\mu$ is central to deformation theory for surfaces.   The Uniformization Theorem provides solutions for the prescribed curvature equation.   The variation theory for the prescribed curvature equation provides an approach for deformations.  The Eells-Sampson variation theory for harmonic maps provides a closely related approach.  Riemann surfaces are also described as identification spaces.  Deformations of geometric structures are correspondingly described as \v{C}ech $1$-cocycles with values in special vector fields.  The \v{C}ech {\em plates} viewpoint provides a context for describing the Fenchel-Nielsen twist deformation, as well as Thurston's earthquake and Penner's shearing.  There is also a description of deformations following the classical theory of automorphic forms; the approach is based on explicit evaluations of integrals and sums - and expressing results in geometric terms, \cite{Wlsymp}.  A basic consideration is a sum over configuration of lines in the hyperbolic plane.  

Many WP results are based on considerations of (complex) differential geometry - for the individual Riemann surfaces - and for Teichm\"{u}ller space. Li-Sun-Yau employ complex differential geometry, analysis and comparison arguments including the Schwarz Lemma of Yau. Recently techniques of synthetic geometry and especially $CAT(0)$ geometry have also become important, \cite{Brkwpvs, BF, BrMs, BMM, DW2, Wlcomp, Wlbhv}.  The intersection theory considerations involve differential forms, symplectic geometry and extensive considerations of Greens functions.  In fact Greens functions appear in the expressions for WP curvature, curvature of the family hyperbolic metric \cite{Wlchern}, formulas of bosonic Polyakov string theory \cite{DP}, and in the form of Eisenstein series in the Takhtajan-Zograf local index formula \cite{TZ}, and in the curvature of the renormalized family hyperbolic metric \cite{Wlcusps}. Understanding Greens functions is basic for understanding curvature.

Topological considerations are also used to investigate WP geometry.  The pants graphs $P(R)$, curve complex $C(R)$ and measured geodesic laminations play important roles.  Mirzakhani's considerations involve Thurston's symplectic form on the space of measured geodesic laminations.  Brock-Masur-Minsky employ methods developed in examination of Thurston's ending lamination conjecture.

Dynamical considerations are also used.  McMullen's approach reveals a connection between Hausdorff dimension, norms of holomorphic forms, and the central limit theorem for geodesic flows, especially the variance of observables of mean zero, \cite{McMthr}.  The Brock-Minsky-Masur considerations on limits of sequences in $C(R)$ include beginning dynamical arguments. 

\section{Lecture material}
\label{lect}

\begin{itemize}
\item Basics: hyperbolic surfaces, deformations and topological methods (2 lectures)

\begin{itemize}
\item Teichm\"{u}ller space, models for infinitesimal deformations; Weil-Petersson and 
Teichm\"{u}ller metrics; mapping class groups; moduli space; Deligne-Mumford compactification
\item constructing surfaces and hyperbolic metrics; enhanced collar lemma and simple geodesics; Fenchel-Nielsen deformation; geodesic-length functions; twist-length duality; first and second twist derivatives of geodesic-length functions
\item measured geodesic laminations; Thurston's symplectic form; curve complex and pants graph; Gromov boundary of the curve complex; ending lamination methods
\end{itemize}

\item Estimating coset sums of inverse-exponential-distance (including estimating Greens functions)
\begin{itemize}
\item examples of inverse-exponential sums; enhanced mean value estimate; estimating on thick and thin regions
\end{itemize} 
\item Geodesic-length functions (1 1/2 lectures)
\begin{itemize}
\item cosine formula; Riera's formula; Hessian of geodesic-length; convexity; small-length expansions; special Schwarz lemma; expansion of WP Levi Civita connection; applications
\end{itemize}
\item Augmented Teichm\"{u}ller space $\Tbar$ and $CAT(0)$ (1 1/2 lectures)
\begin{itemize}
\item definition $\Tbar$; $CAT(0)$; Brock's approximation of geodesics; Brock's quasi isometry; $\Tbar$ as an infinite polyhedra; $MCG\,=\,\operatorname{Isom}_{WP}$; classification of flats; Alexandrov tangent cones
\end{itemize}
\item Mirzakhani's volume recursion and {\em prime geodesic theorem} (1 1/2 lectures)
\begin{itemize}
\item Mirzakhani-McShane identity; the recursion of volume integrals; simple geodesics and counting lattice points on MGL; Makover-McGowan's expected short length statistics  
\end{itemize}
\item Behavior of WP geodesics and topological transitivity of geodesic flow
\begin{itemize}
\item compactness theorem for geodesics; ending laminations for infinite WP geodesics;  recurrent WP geodesics; encoding geodesics
\end{itemize}
\item Curvature and index formulas for families hyperbolic metrics (1 1/2 lectures)
\begin{itemize}
\item formulas; perturbations of the prescribed curvature equation; degenerating curvatures; the plumbing family $\{zw=t\} \rightarrow \{t\}$; families of holomorphic $2$-differentials; Mumford's good \& pre-log-log metrics 

\end{itemize}
\end{itemize}
\pagebreak
\nocite{Wlext}
\nocite{MaMc}

\begin{thebibliography}{PWW10}

\bibitem[Ahl61]{Ahsome}
Lars~V. Ahlfors.
\newblock Some remarks on {T}eichm\"uller's space of {R}iemann surfaces.
\newblock {\em Ann. of Math. (2)}, 74:171--191, 1961.

\bibitem[BF06]{BF}
Jeffrey Brock and Benson Farb.
\newblock Curvature and rank of {T}eichm\"uller space.
\newblock {\em Amer. J. Math.}, 128(1):1--22, 2006.

\bibitem[BH99]{BH}
Martin~R. Bridson and Andr{\'e} Haefliger.
\newblock {\em Metric spaces of non-positive curvature}.
\newblock Springer-Verlag, Berlin, 1999.

\bibitem[BM07]{BrMs}
Jeffrey Brock and Howard Masur.
\newblock Coarse and synthetic {W}eil-{P}etersson geometry: quasi-flats,
  geodesics, and relative hyperbolicity.
\newblock preprint, 2007.

\bibitem[BMM07]{BMM}
Jeffrey Brock, Howard Masur, and Yair Minsky.
\newblock Asymptotics of {W}eil-{P}etersson geodesics {I}: ending laminations,
  recurrence, and flows.
\newblock preprint, 2007.

\bibitem[Bro03]{Brkwp}
Jeffrey~F. Brock.
\newblock The {W}eil-{P}etersson metric and volumes of 3-dimensional hyperbolic
  convex cores.
\newblock {\em J. Amer. Math. Soc.}, 16(3):495--535 (electronic), 2003.

\bibitem[Bro05]{Brkwpvs}
Jeffrey~F. Brock.
\newblock The {W}eil-{P}etersson visual sphere.
\newblock {\em Geom. Dedicata}, 115:1--18, 2005.

\bibitem[DP86]{DP}
Eric D'Hoker and D.~H. Phong.
\newblock Multiloop amplitudes for the bosonic {P}olyakov string.
\newblock {\em Nuclear Phys. B}, 269(1):205--234, 1986.

\bibitem[DW03]{DW2}
Georgios Daskalopoulos and Richard Wentworth.
\newblock Classification of {W}eil-{P}etersson isometries.
\newblock {\em Amer. J. Math.}, 125(4):941--975, 2003.

\bibitem[FiM09]{GFM}
G{\'e}rard Freixas~i Montplet.
\newblock An arithmetic {R}iemann-{R}och theorem for pointed stable curves.
\newblock {\em Ann. Sci. \'Ec. Norm. Sup\'er. (4)}, 42(2):335--369, 2009.

\bibitem[GGHar]{Zh4}
Ren Guo, Subhojoy Gupta, and Zheng Huang.
\newblock Curvatures on the {T}eichm\"{u}ller curve.
\newblock {\em Indiana Univ. Math. Jour.}, to appear.

\bibitem[Gol86]{Gdsymp}
William~M. Goldman.
\newblock Invariant functions on {L}ie groups and {H}amiltonian flows of
  surface group representations.
\newblock {\em Invent. Math.}, 85(2):263--302, 1986.

\bibitem[Hua05]{Zh}
Zheng Huang.
\newblock Asymptotic flatness of the {W}eil-{P}etersson metric on
  {T}eichm\"uller space.
\newblock {\em Geom. Dedicata}, 110:81--102, 2005.

\bibitem[Hua07a]{Zh2}
Zheng Huang.
\newblock On asymptotic {W}eil-{P}etersson geometry of {T}eichm\"uller space of
  {R}iemann surfaces.
\newblock {\em Asian J. Math.}, 11(3):459--484, 2007.

\bibitem[Hua07b]{Zh3}
Zheng Huang.
\newblock The {W}eil-{P}etersson geometry on the thick part of the moduli space
  of {R}iemann surfaces.
\newblock {\em Proc. Amer. Math. Soc.}, 135(10):3309--3316 (electronic), 2007.

\bibitem[LSY04]{LSY1}
Kefeng Liu, Xiaofeng Sun, and Shing-Tung Yau.
\newblock Canonical metrics on the moduli space of {R}iemann surfaces. {I}.
\newblock {\em J. Differential Geom.}, 68(3):571--637, 2004.

\bibitem[LSY05]{LSY2}
Kefeng Liu, Xiaofeng Sun, and Shing-Tung Yau.
\newblock Canonical metrics on the moduli space of {R}iemann surfaces. {II}.
\newblock {\em J. Differential Geom.}, 69(1):163--216, 2005.

\bibitem[McM06]{McMthr}
Curtis~T. McMullen.
\newblock Thermodynamics, dimension and the {W}eil-{P}etersson metric.
\newblock preprint, 2006.

\bibitem[Mir07a]{Mirvol}
Maryam Mirzakhani.
\newblock Simple geodesics and {W}eil-{P}etersson volumes of moduli spaces of
  bordered {R}iemann surfaces.
\newblock {\em Invent. Math.}, 167(1):179--222, 2007.

\bibitem[Mir07b]{Mirwitt}
Maryam Mirzakhani.
\newblock Weil-{P}etersson volumes and intersection theory on the moduli space
  of curves.
\newblock {\em J. Amer. Math. Soc.}, 20(1):1--23 (electronic), 2007.

\bibitem[Mir08]{Mirgrow}
Maryam Mirzakhani.
\newblock Growth of the number of simple closed geodesics on hyperbolic
  surfaces.
\newblock {\em Ann. of Math. (2)}, 168(1):97--125, 2008.

\bibitem[MM05]{MaMc}
Eran Makover and Jeffrey McGowan.
\newblock The length of closed geodesics on random {R}iemann surfaces.
\newblock Arxiv:math/0504175, 2005.

\bibitem[Pen93]{Pen2frm}
Robert~C. Penner.
\newblock The {P}oincar\'e dual of the {W}eil-{P}etersson {K}\"ahler two-form.
\newblock {\em Comm. Anal. Geom.}, 1(1):43--69, 1993.

\bibitem[PWW10]{PWWl}
Mark Pollicott, Howard Weiss, and Scott~A. Wolpert.
\newblock Topological dynamics of the {W}eil-{P}etersson geodesic flow.
\newblock {\em Adv. Math.}, 223(4):1225--1235, 2010.

\bibitem[Roy75]{Royicm}
H.~L. Royden.
\newblock Intrinsic metrics on {T}eichm\"uller space.
\newblock In {\em Proceedings of the International Congress of Mathematicians
  (Vancouver, B. C., 1974), Vol. 2}, pages 217--221. Canad. Math. Congress,
  Montreal, Que., 1975.

\bibitem[Tro86]{Trcurv}
A.~J. Tromba.
\newblock On a natural algebraic affine connection on the space of almost
  complex structures and the curvature of {T}eichm\"uller space with respect to
  its {W}eil-{P}etersson metric.
\newblock {\em Manuscripta Math.}, 56(4):475--497, 1986.

\bibitem[TZ87]{TZcomp}
L.~A. Takhtadzhyan and P.~G. Zograf.
\newblock A local index theorem for families of {$\overline\partial$}-operators
  on {R}iemann surfaces.
\newblock {\em Uspekhi Mat. Nauk}, 42(6(258)):133--150, 248, 1987.

\bibitem[TZ88]{TZpunc}
L.~A. Takhtajan and P.~G. Zograf.
\newblock The {S}elberg zeta function and a new {K}\"ahler metric on the moduli
  space of punctured {R}iemann surfaces.
\newblock {\em J. Geom. Phys.}, 5(4):551--570 (1989), 1988.

\bibitem[TZ91]{TZ}
L.~A. Takhtajan and P.~G. Zograf.
\newblock A local index theorem for families of $\overline\partial$-operators
  on punctured {R}iemann surfaces and a new {K}\"ahler metric on their moduli
  spaces.
\newblock {\em Comm. Math. Phys.}, 137(2):399--426, 1991.

\bibitem[\Wl89]{Wfharm}
Michael \WlName{Wolf}.
\newblock The {T}eichm\"uller theory of harmonic maps.
\newblock {\em J. Differential Geom.}, 29(2):449--479, 1989.

\bibitem[\Wp83a]{Wlhomol}
Scott~A. \WpName{Wolpert}.
\newblock On the homology of the moduli space of stable curves.
\newblock {\em Ann. of Math. (2)}, 118(3):491--523, 1983.

\bibitem[\Wp83b]{Wlsymp}
Scott~A. \WpName{Wolpert}.
\newblock On the symplectic geometry of deformations of a hyperbolic surface.
\newblock {\em Ann. of Math. (2)}, 117(2):207--234, 1983.

\bibitem[\Wp85a]{Wlpos}
Scott~A. \WpName{Wolpert}.
\newblock On obtaining a positive line bundle from the {W}eil-{P}etersson
  class.
\newblock {\em Amer. J. Math.}, 107(6):1485--1507 (1986), 1985.

\bibitem[\Wp85b]{Wldtau}
Scott~A. \WpName{Wolpert}.
\newblock On the {W}eil-{P}etersson geometry of the moduli space of curves.
\newblock {\em Amer. J. Math.}, 107(4):969--997, 1985.

\bibitem[\Wp86a]{Wlchern}
Scott~A. \WpName{Wolpert}.
\newblock Chern forms and the {R}iemann tensor for the moduli space of curves.
\newblock {\em Invent. Math.}, 85(1):119--145, 1986.

\bibitem[\Wp86b]{Wlthur}
Scott~A. \WpName{Wolpert}.
\newblock Thurston's {R}iemannian metric for {T}eichm\"uller space.
\newblock {\em J. Differential Geom.}, 23(2):143--174, 1986.

\bibitem[\Wp87]{Wlnielsen}
Scott~A. \WpName{Wolpert}.
\newblock Geodesic length functions and the {N}ielsen problem.
\newblock {\em J. Differential Geom.}, 25(2):275--296, 1987.

\bibitem[\Wp03]{Wlcomp}
Scott~A. \WpName{Wolpert}.
\newblock {G}eometry of the {W}eil-{P}etersson completion of {T}eichm\"{u}ller
  space.
\newblock In {\em Surveys in Differential Geometry VIII: Papers in Honor of
  Calabi, Lawson, Siu and Uhlenbeck}, pages 357--393. Intl. Press, Cambridge,
  MA, 2003.

\bibitem[\Wp07]{Wlcusps}
Scott~A. \WpName{Wolpert}.
\newblock Cusps and the family hyperbolic metric.
\newblock {\em Duke Math. J.}, 138(3):423--443, 2007.

\bibitem[\Wp08]{Wlbhv}
Scott~A. \WpName{Wolpert}.
\newblock Behavior of geodesic-length functions on {T}eichm\"uller space.
\newblock {\em J. Differential Geom.}, 79(2):277--334, 2008.

\bibitem[\Wp09]{Wlext}
Scott~A. \WpName{Wolpert}.
\newblock Extension of the {W}eil-{P}etersson connection.
\newblock {\em Duke Math. J.}, 146(2):281--303, 2009.

\end{thebibliography}

\providecommand\WlName[1]{#1}\providecommand\WpName[1]{#1}\providecommand\Wl{W%
lf}\providecommand\Wp{Wlp}\def\cprime{$'$}

\end{document}